\documentclass{amsart}
\usepackage{a4}

\newtheorem{theorem}{Theorem}[section]

\newtheorem{proposition}[theorem]{Proposition}
\newtheorem{corollary}[theorem]{Corollary}

\theoremstyle{definition}

\theoremstyle{remark}

\numberwithin{equation}{section}

\def\xceil[#1]{{{\left\ulcorner{#1}\right\urcorner}}}
\def\floor[#1]{{{\left\lfloor{#1}\right\rfloor}}}
\def\ceil[#1]{{{\left\lceil{#1}\right\rceil}}}
\def\brace#1{{{\left\lbrace{#1}\right\rbrace}}}
\def\paren(#1){{{\left({#1}\right)}}}
\def\Z{{\mathbb Z}}
\def\congruent{\equiv}
\DeclareMathOperator*{\lcm}{lcm}

\newcommand{\therosteritem}[1]{{(#1)}}

\begin{document}

\title{Economical numbers}

\author{R.G.E. Pinch}
\address{Queens' College, Silver Street,
Cambridge  CB3 9ET, U.K.}
\email{rgep@cam.ac.uk}


\date{9 February 1998}


\begin{abstract}
A number $n$ is said to be {\it economical} if the prime power factorisation of $n$
can be written with no more digits than $n$ itself.  We show that under a plausible
hypothesis, related to the twin prime conjecture, there are arbitrarily long sequences
of consecutive economial numbers, and exhibit such a sequence of length 9.
\end{abstract}

\maketitle

\section{Introduction}
In \cite{San:equidigit}, Bernardo Recam\'an Santos defined a number $n$ to be
{\it equidigital}\/\footnote{We prefer this spelling to ``equadigital''.}
if the prime power factorisation of $n$ requires the same number of decimal
digits as $n$, and {\it economical} if the prime power factorisation requires
no more digits.  He asked whether there were arbitrarily long
sequences of consecutive economical numbers.  In \cite{Hes:equidigit}, 
Richard Hess observed
that the five consecutive integers $13$, $14 = 2 \cdot 7$, $15 = 3 \cdot 5$,
$16 = 2^4$ and $17$ are economical.  We note that the sequence $157$, $158 = 2 \cdot 79$,
$159 = 3 \cdot 53$, $160 = 2^5 \cdot 5$, $161 = 7 \cdot 23$, $162 = 2 \cdot 3^4$,
$163$ is of length seven.  Hess also suggested that there was a longest such string.

In this article, we
shall show that there are strings of consecutive economical numbers of arbitrary
length if Dickson's conjecture, on the simultaneous primality of linear forms,
is true.  This is a generous assumption, since the conjecture includes the twin
prime conjecture as a special case.  We use the idea of the proof to find a
sequence of nine consecutive economical numbers from  1034429177995381247
to 10344291779953812455.

\section{Economical and frugal numbers}

There is no need to restrict to decimal expansions, so we begin by working with
a general base $B$.

Let $\delta(n)$ denote the number of digits required to write $n$ in base $B$,
so that $\delta(n) = k$ if and only if $B^{k-1} \le n < B^k$.  We can express
this conveniently as $\delta(n) = \xceil[\log_B n]$, 
where\footnote{Note that this differs from the usual $\ceil[x]$ when $x$ is an 
integer: indeed, $\xceil[x] = \floor[x]+1$.}
$\xceil[x] = \min\brace{ k \in \Z : k > x }$.

We let $\phi(n)$
denote the number of digits required to write down the prime power factorisation
of $n$, so that if $n = \prod_{i=1}^r p_i^{a_i}$ 
then $\phi(n) = \sum_{i=1}^r \delta(p_i) + \delta'(a_i)$, where $\delta'(a) = \delta(a)$
for $a > 1$ and $\delta'(1) = 0$.  We let $h(n) = \delta(n) - \phi(n)$.

We define $n$ to be {\it equidigital (base $B$)} if $h(n) = 0$,
{\it economical (base $B$)} if $h(n) \ge 0$ and {\it frugal 
(base $B$)} if $h(n) > 0$.  For $k > 0$,
if $h(n) \ge k$ we say that $n$ is {\it $k$--frugal}, so that 1-frugal is the same as
frugal.  We include $1$ as a frugal number.

We begin with some obvious properties.

\begin{proposition}\label{prop:1}
\begin{enumerate}
\item
A prime number is equidigital in every base;
\item
$\max\brace{\delta(m), \delta(n)} \le \delta(m+n) 
           \le 1 + \max\brace{\delta(m), \delta(n)}$;
\item
$\delta(m) + \delta(n) - 1 \le \delta(mn) \le \delta(m) + \delta(n)$;
\item
$b \delta(n) - b \le \delta\paren(n^b) \le b \delta(n)$;
\item
If $m$ and $n$ are coprime, then $\phi(mn) = \phi(m) + \phi(n)$;
\item
$\phi(mn) \le \phi(m) + \phi(n)$;
\item
$\phi(n^b) \le \phi(n) + \delta(b) \log_2 n$;
\item
$h(mn) \le h(m) + h(n) - 1$;
\item
If $m$ is frugal and $n$ is economical, then $mn$ is economical;
\item
If $m$ is frugal then $h(mn) \le h(n)$.
\end{enumerate}
\end{proposition}
\begin{proof}
Part \therosteritem1 follows immediately from the definition.

For part \therosteritem2, suppose that $m \le n$ and $\delta(n) = k$.  Then
$B^{k-1} < m+n \le 2n \le 2 B^k \le B^{k+1}$.

For part \therosteritem3, we note that if $\delta(m) = l$ and $\delta(n) = k$
then $B^{l-1} \le m < B^l$ and $B^{k-1} \le n < B^k$, so 
that $B^{k+l-2} \le mn < B^{k+l}$.  Part \therosteritem4 follows by induction.

For part \therosteritem5, we note that if $m$ and $n$ are coprime, then their
prime power factorisation are disjoint.

For part \therosteritem6, we need to consider only the prime power factors
that $m$ and $n$ have in common.  If $p^a$ occurs in $m$ and $p^b$ in $n$,
then $p^{a+b}$ occurs in $mn$.  Assume that $a \ge b$.
We claim that $\delta'(a+b) + \delta(p) \le \delta'(a) + \delta'(b) + 2\delta(p)$,
that is, $\delta'(a+b) \le \delta'(a) + \delta'(b) + \delta(p)$.  
If $a = b = 1$, then $\delta'(a+b) = \delta'(2) = \delta(2) \le \delta(p)$.
If $a > b = 1$, then $\delta'(a+b) = \delta'(a+1) = \delta(a+1) \le \delta(a) + 1$ 
by \therosteritem2.
If $a \ge b > 1$ then $\delta'(a+b) = \delta(a+b) \le \delta(a) + \delta(b)$,
again by \therosteritem2.  

For part \therosteritem7, let $n = \prod_{i=1}^r p_i^{a_i}$.  Assume $b > 1$.
Then $\phi(n^b) = \sum_{i=1}^r \delta\paren(p_i) + \delta\paren(a_i b)$ and using
\therosteritem3 this is at 
most $\sum_{i=1}^r \delta\paren(p_i) + \delta(b) +  \delta'\paren(a_i) 
= \phi(n) + r \delta(b)$.
Now $r \le \log_2 n$, so $\phi(n^b) \le \phi(n) + \delta(b) \log_2 n$.

We obtain \therosteritem8 by combining \therosteritem 3 and \therosteritem6:
parts \therosteritem9 and \therosteritem{10} are special cases of \therosteritem8.
\end{proof}

We use the following proposition to prove the result we shall use in the next section.

\begin{proposition}\label{prop:2}
Fix a base $B$.  
\begin{enumerate}
\item
Suppose that $r$ and $s$ are coprime.  For each $k \ge 0$
there are infinitely many $k$-frugal numbers coprime to $s$ and divisible
exactly by the prime powers in $r$.
\item
Given $r$, for each $k \ge 0$ there are infinitely many $k$-frugal numbers 
divisible by $r$ and containing only the same prime factors as $r$.
\end{enumerate}
\end{proposition}
\begin{proof}
For \therosteritem1, take $p$ to be a prime greater than $rs$ and let $n = r p^a$.  
Assume $a \ge 2$.  We
have $\phi(n) = \phi(r) + \delta(a) + \delta(p)$
and $\delta(n) \ge \delta(r) + \delta\paren(p^a) \ge \delta(r) + a(\delta(p) - 1)$.
So 
\begin{eqnarray*}
h(n) &\ge& \delta(r) + a(\delta(p) - 1) - \phi(r) - \delta(a) - \delta(p)	\\
     & &\quad\quad= a(\delta(p) - 1) - \delta(a) +h(r) - \delta(p)
\end{eqnarray*}
and this tends to infinity with $a$.

For \therosteritem2, consider $n = r^a$ for some $a > 1$.  By
Proposition \ref{prop:1} \therosteritem4 and \therosteritem7 we 
have 
\begin{eqnarray*}
h(n) &=& \delta(r^a) - \phi(r^a) \ge  a \delta(r) - a  
		- \paren({\phi(r) + \delta(a) \log_2 r})	\\
     &\ge& a(\delta(r)-1)  - \delta(a) \log_2 r - \phi(r),
\end{eqnarray*}
and again this tends to infinity with $a$.
\end{proof}

Although we shall not need it, we can prove the converse of 
Proposition \ref{prop:1} \therosteritem{10}.  We use the following
result of Baker and Harman \cite{BH:primegaps}.

\begin{proposition}\label{prop:3}
For  $x$ sufficiently large, there is always a prime between $x$ and $x + x^\theta$, 
where we can take $\theta = 0.535$. \qed
\end{proposition}

For our purposes any $\theta < 1$ would suffice.  See
Ribenboim \cite{Rib:rec},4.II.C for other results in this direction.

\begin{proposition}\label{prop:4}
Suppose that $m$ has the property that $mn$ is economical whenever $n$ is.
Then $m$ is frugal.
\end{proposition}
\begin{proof}
Suppose that $m$ has the property stated, and let $\delta(m) = l$.  Then
$B^{l-1} \le m < B^l$, so $\frac{B^l}{m} > 1$.  Let $k$ be large enough that 
$(B^k + B^{\theta k}) / B^k = 1 + B^{(\theta-1)k} < B^l / m$: we may also
assume that $k > l$ and that $B^k$ is large enough for Proposition \ref{prop:3}
to hold.  Take a prime $p$ in the interval $B^k < p \le B^k + B^{\theta k}$: 
by Proposition \ref{prop:1} \therosteritem1, $p$ is economical, 
so $mp$ is economical by the assumption on $m$.  
But $mp < m\paren({1 + B^{(\theta-1)k}}) < B^{l+k}$, 
so $\delta(mp) = \delta(m) + \delta(p) - 1 = k+l$.  	
We have $p > m$, so $p$ is coprime to $m$ and 
hence $\phi(mp) = \phi(m) + \phi(p) = \phi(m) + \delta(p) = \phi(m) + k+1$.
So $h(mp) = \delta(mp) - \phi(mp) = k + l - \phi(m) - k - 1 = l - \phi(m) - 1 = h(m) - 1$.
But by our assumptions, $mp$ is economical, so that $h(mp) \ge 0$ and so we must
have $h(m) \ge 1$: that is, $m$ is frugal.
\end{proof}

\section{Extravagant numbers}

We can use Proposition \ref{prop:3} to show that there are numbers $n$ 
with $h(n) \ll 0$: presumably we should call these {\it extravagant}.

\begin{proposition}\label{prop:5}
Fix a base $B$.  For $k > 0$ there are infinitely many numbers with $h(n) \le -k$.
\end{proposition}
\begin{proof}
Choose $t$ large enough that $B^t$ exceeds the bound of Proposition \ref{prop:3}, and also
$B^{t(1-\theta)} > 2^\theta (k+2)^2$.  Put $L = 1 + \frac{2^\theta}{B^{t(1-\theta)}}$.
Then 
$$
L  < 1 + \frac{1}{(k+2)^2}
	< 1 + \frac{1}{2(k+1)}
$$
so
$$
L^{k+1} < \paren(1 + \frac{1}{2(k+1)})^{k+1}
	< \exp\paren({1/2}) < 2 ,
$$
since $\paren(1 + x/n)^n$ tends to $e^x$ from below for $x > 0$.    

Consider the intervals $[B^t L^i, B^t L^{i+1}]$, for $i = 0,\ldots,k$.
The length of the $i$-th interval is $B^t L^i(L-1)$ and 
$$
B^t L^i (L-1) \ge B^t (L-1) = 2^\theta B^{t \theta}
	> L^{k+1} B^{t \theta} \ge L^i B^{t \theta} > (L^i B^t)^\theta ,
$$
so that Proposition \ref{prop:3} applies and each interval contains a prime, say $p_i$.
Further, $B^k L^{k+1} < 2 B^k \le B^{k+1}$.  Hence $\delta(p_i) = t+1$ for each $i$.

If we put $n = \prod_{i=0}^k p_i$, 
then $\phi(n) = \sum_{i=0}^k \delta(p_i) = (k+1)(t+1)$.
But 
\begin{eqnarray*}
n &<& B^{t(k+1)} L^{(k+1)(k+2)/2} 
                    < B^{t(k+1)} \paren(1 + \frac{1}{(k+1)(k+2)})^{(k+1)(k+2)/2} \\
		&<& B^{t(k+1)} \exp\paren({1/2}) < B^{t(k+1)+1},
\end{eqnarray*}
so that $\delta(n) = t(k+1)+1$.  So $h(n) = t(k+1) + 1 - (k+1)(t+1) = -k$.
\end{proof}

We note that the idea of the proof can be applied to intervals of the form
$[B^t L^{-i-1}, B^t L^{-i}]$ to show that there are infinitely many squarefree
economical numbers with distinct prime factors, providing an alternative solution
to \cite{San:equidigit} in any base.

\section{Dickson's conjecture}

Dickson \cite{Dic:conjecture} conjectured that a family of linear functions
$f_i(n) = a_i n + b_i$, $i = 1,...,t$, with integer coefficients $a_i$, $b_i$,
would be simultaneously prime
unless they ``obviously cannot'': that is, unless there is an integer $m > 1$
such that $m$ divides the product $f_1(n)f_2(n)\cdots f_t(n)$ for every value 
of $n$.

This rather powerful conjecture would imply the {\it twin primes conjecture},
that there are infinitely many twin prime pairs $(p,p+2)$, and the 
{\it Sophie Germain primes conjecture}, that there are infinitely primes $p$ for
which $2 p+1$ is also 
prime\footnote{If $p$ is a Sophie Germain prime, then $2p+1$ is a {\it safe} prime.}.

Schinzel \cite{SS:hypH} extended Dickson's conjecture to the analogous 
``Hypothesis H'' for arbitrary integer polynomials $f_i(x)$.

\section{Consecutive economical numbers}

\begin{theorem}
Fix a base $B$, and $k \ge 0$.
If Dickson's conjecture holds, there are arbitrarily long sequences of 
consecutive $k$-frugal numbers.
\end{theorem}
\begin{proof}
Suppose we wish to find a sequence of $t$ consecutive $k$-frugal numbers,
say $N, N+1, \ldots, N+t-1$.  Let $p_h$ be the largest prime less than $t$
and for the primes $p_i \le p_h$ let $p_i^{a_i}$ be the largest power 
of $p_i$ less than $t$.

We shall insist that $N$ be divisible by $f_0 = \prod_{i=1}^h p_i^{a_i + 1}$.
This implies that the $N+j$ must be divisible by certain powers of the $p_i$:
in fact, for $j=1,\ldots,t-1$, if $p_i^{b_{i,j}}$ exactly divides $j$ then the 
same power of $p_i$ must exactly divide $N+j$.  Let $f_j = \prod_i p_i^{b_{i,j}}$.
(In fact, $f_j = j$ but we shall not use that.)  We have $f_1 = 1$.
For each $j = 2,\ldots,t-1$ we use Proposition \ref{prop:2} \therosteritem1 to 
find $m_j$ such that
$h(m_j) \ge k+2 - h(f_j) \ge 1$ and each $m_j$ is coprime to $f_0$ and the
other $m_{j'}$.  We thus ensure that $m_j f_j$ is $(k+1)$-frugal.

Since $f_0$ and the $m_j$ are coprime, we can use the Chinese Remainder Theorem
to solve the simultaneous
congruences $N+j \congruent 0 \bmod f_j m_j$ to find infinitely many
solutions of the form $N \congruent N_0 \bmod M$ 
where $M = \lcm\brace{f_0, m_j}$.  The cofactors $C_j = \frac{N+j}{f_j m_j}$
will satisfy congruences of the 
form $C_j \congruent S_j \bmod M_j$ where $S_j = \frac{N_0+j}{f_j m_j}$ 
and $M_j = \frac{M}{f_j m_j}$.

By Dickson's conjecture, the forms $C_j = S_j + M_j x$ will be simultaneously prime 
for infinitely many values for $x$.  For such values, we 
have $h\paren(f_j m_j) \ge k+1$
and $C_j$ economical, so by Proposition \ref{prop:1} \therosteritem8 
the values $N+j = f_j m_j C_j$ are all $k$-frugal,
as required.
\end{proof}

\begin{corollary}
Fix a base $B$.
If Dickson's conjecture holds, there are arbitrarily long sequences of 
consecutive economical numbers.\qed
\end{corollary}

In practice there are modifications one could make to the construction
given in the proof of the Theorem.  For example, it may be more
efficient to take the sequence of numbers $N-1,\ldots,N+t-2$, since none of
the $p_i$ will divide $N-1$ and we can take $f_{-1} = 1$.

A further variant would be to take $m_0$ to be composed of further 
powers of the primes $p_i$, using Proposition \ref{prop:2} \therosteritem2 rather than
\therosteritem1 to ensure that $f_0 m_0$ is frugal: for example, we might
take $m_0$ to be a power of $f_0$.

\section{Examples}

We illustrate with some example.  We take base $B = 10$.
Suppose we wish to find seven consecutive economical
numbers.  We have $p_3 = 5$ and $2^2, 3^1, 5^1$ the highest powers less than $7$,
so that $f_0 = 2^3 \cdot 3^2 \cdot 5^2 = 1800$.  
We find that $f_j = j$, $j=1,\ldots,6$.  We take $m_0 = 7^6$, 
$m_1 = 1$, $m_2 = 11^4$, $m_3 = 13^4$,
$m_4 = 17^4$, $m_5 = 19^3$, $m_6 = 23^4$.  
We have     $M = 14196220211350791776356766371800$
and $N \congruent N_0 = 5599355285926686611723646146400 \bmod M$.

For $N = 133365337188083812598934543492600 = N_0 + 9M$ we have
{\small
\begin{eqnarray*}
133365337188083812598934543492599 &=& 29 \cdot 598426817561 \cdot 7684823934473500571\\
133365337188083812598934543492600 &=& 2^{3} \cdot 3^{2} \cdot 5^{2} \cdot 7^{6} \cdot 649567 \cdot 969523340521703729\\
133365337188083812598934543492601 &=& 133365337188083812598934543492601\\
133365337188083812598934543492602 &=& 2 \cdot 11^{4} \cdot 89 \cdot 113 \cdot 6701 \cdot 67582485977340653273\\
133365337188083812598934543492603 &=& 3 \cdot 13^{4} \cdot 170173187 \cdot 832591957 \cdot 10985629799\\
133365337188083812598934543492604 &=& 2^{2} \cdot 17^{4} \cdot 73 \cdot 16141 \cdot 338792660124002968867\\
133365337188083812598934543492605 &=& 5 \cdot 19^{3} \cdot 46399 \cdot 188563981 \cdot 444472402866601\\
133365337188083812598934543492606 &=& 2 \cdot 3 \cdot 23^{5} \cdot 3453444945058703174533907
\end{eqnarray*}
}

We find that $N-1,\ldots,N+6$ form a sequence of eight consecutive economical
numbers: even though we only wanted seven! We note that not all the cofactors 
are prime, although they are all economical.

If we adopt the variant of choosing $m_0$ to be a power of $f_0$, then in this
example we would have $m_0 = f_0 = 1800$, $m_2 = 7^4$, $m_3 = 11^4$, $m_4 = 13^4$,
$m_5 = 17^3$ and $m_6 = 19^4$.  We have $M = 2082775632877914851396520000$ 
and $N \congruent N_0 = 1625787524296851742054440000 \bmod M$.

For $N = 1625787524296851742054440000$ we have the sequence of length seven
\begin{eqnarray*}
1625787524296851742054440000&=&2^{6} \cdot 3^{4} \cdot 5^{4} \cdot 501786272931127080881 \\
1625787524296851742054440001&=&2237 \cdot 726771356413433948169173 \\
1625787524296851742054440002&=&2 \cdot 7^{4} \cdot 678077 \cdot 499301208591113813 \\
1625787524296851742054440003&=&3 \cdot 11^{4} \cdot 23 \cdot 37030463 \cdot 43459508766889 \\
1625787524296851742054440004&=&2^{2} \cdot 13^{4} \cdot 183907 \cdot 77380605917154763 \\
1625787524296851742054440005&=&5 \cdot 17^{3} \cdot 227 \cdot 65239 \cdot 4469036127622909 \\
1625787524296851742054440006&=&2 \cdot 3 \cdot 19^{4} \cdot 53 \cdot 6297281 \cdot 1751557 \cdot 3556681
\end{eqnarray*}

Again we note that not all the cofactors are prime, although they are all
economical.

The precise assignment of the powers of the small primes $p_i$ to the number $N+j$
in the proof of the Theorem is quite arbitrary.  If we have different elements of
the sequence divisible by the highest powers of each $p_i$, it seems that we can
do better: 

\begin{eqnarray*}
1034429177995381247  &=&  51551 \cdot 20066132140897 			\\
1034429177995381248  &=&  2^{9} \cdot 3 \cdot 88651 \cdot 7596716293 	\\
1034429177995381249  &=&  17 \cdot 60848775176198897 			\\
1034429177995381250  &=&  2 \cdot 5^{5} \cdot 76757 \cdot 2156268073 	\\
1034429177995381251  &=&  3^{5} \cdot 19 \cdot 383 \cdot 584981475541 	\\
1034429177995381252  &=&  2^{2} \cdot 7^{5} \cdot 154267 \cdot 99741877 	\\
1034429177995381253  &=&  394007 \cdot 2625408122179 			\\
1034429177995381254  &=&  2 \cdot 3 \cdot 23 \cdot 59^{3} \cdot 1367 \cdot 26699131 	\\
1034429177995381255  &=&  5 \cdot 26264543 \cdot 7877001157 .
\end{eqnarray*}

\section{Distribution of frugal numbers}

Again we take base $B = 10$.  We computed $h(n)$ for $2 \le n \le 1000000$.  

The number of integers up to $5 \cdot 10^8$ which have $h(n) = k$ is
{\small
$$\begin{array}{rrrrrrrrrrrrr}
k =-6 & -5     &   -4    &    -3    &     -2    &     -1    \\
 1313 & 195341 & 5101112 & 44435592 & 153988692 & 208380123 \\
\end{array}
$$
$$\begin{array}{rrrrrrrrrrrrr}
k =  0   &    1    &   2    &  3    & 4    & 5   & 6  \\
86441875 & 1297001 & 140575 & 16670 & 1483 & 207 & 16 \\
\end{array}
$$
}

The smallest 6-frugal number is $40353607 = 7^9$: the smallest 6-extravagant
number is $8314020 = 2^2 3^2 5 \cdot 11 \cdot 13 \cdot 17 \cdot 19$.

We examined consecutive economical numbers up to $10^6$.
The longest strings of consecutive economical numbers we found were of
length 7, starting at 157, 108749, 109997, 121981, 143421.  The longest
strings of consecutive frugal numbers were only of length two, the first 
starting with $4374$.

{
}
\ifx\undefined\bysame
\newcommand{\bysame}{\leavevmode\hbox to3em{\hrulefill}\,}
\fi

\end{document}